\documentclass[12pt, parskip=full]{article}
\usepackage{amsmath}
\usepackage{amsfonts}
\usepackage{amssymb}
\usepackage{graphicx}
\renewcommand{\bf}{\bfseries}
\renewcommand{\it}{\itshape}
\def\a{\alpha}

\def\be{\begin{equation}}
\def\ee{\end{equation}}

\newtheorem{Theorem}{Theorem}[section]

\newtheorem{Lemma}[Theorem]{Lemma}

\newtheorem{Definition}[Theorem]{Definition}

\def\proof#1. {\par
                      \ifdim\lastskip<15pt
                      \removelastskip\penalty-200
                      \vskip5pt plus3pt minus3pt
                      \fi
                       {\def\a{#1}
                       \ifx\a\empty
                       {\noindent\bf Proof.}
                       \else
                       {\noindent\bf Proof of #1.}
                       \fi}\enspace}

\def\endproof{\hfill\hspace{-6pt}\rule[-4pt]{6pt}{6pt}
\vskip8pt plus3pt minus 3pt}

\allowdisplaybreaks[4]

\begin{document}
\raggedright
\title{Bounds for Extreme Zeros of Quasi-orthogonal Ultraspherical Polynomials}

\author{Kathy Driver and Martin E. Muldoon}

\date{February 9, 2016}
\maketitle
\begin{abstract}

We discuss and compare upper and lower bounds obtained by two different methods
for the positive zero of the ultraspherical polynomial $C_{n}^{(\lambda)}$ that is
greater than $1$ when $-3/2 < \lambda < -1/2.$ Our first approach uses mixed three term
recurrence relations and interlacing of zeros while the second approach uses a method 
going back to Euler and Rayleigh and already applied to Bessel functions and Laguerre and $q$-Laguerre polynomials. We use the bounds obtained by the second method
to simplify the proof of the interlacing  of the zeros
of $(1-x^2)C_{n}^{(\lambda)}$ and $C_{n+1}^{(\lambda)}$, for $-3/2 < \lambda < \infty$.

\end{abstract}





\section{Introduction}
For $\lambda > -1/2$, the sequence of ultraspherical polynomials $\{C_{n}^{(\lambda)}\} _{n=0}^\infty$ is
orthogonal on $[-1,1]$ with respect to the positive measure $(1-x^2)^{\lambda -1/2}$ and all
the zeros of  $C_{n}^{(\lambda)}$ lie in $(-1,1).$  As the parameter $\lambda$ decreases through
$-1/2,$ the zeros  of $C_{n}^{(\lambda)}$ depart from the interval $(-1,1)$ in pairs through the endpoints
as $\lambda$ decreases through the values $-1/2, -3/2,\dots,
-\lfloor n/2\rfloor +1/2$.  Here, we pay particular attention to the case  $-3/2 < \lambda < -1/2,$ where $\lambda$ is fixed and
it is known \cite[Cor. 2]{BDR} that the positive zeros, listed in decreasing order, satisfy
\begin{equation}
x_{\lfloor n/2 \rfloor,n}(\lambda) < \cdots < x_{2,n}(\lambda) < 1 < x_{1,n}(\lambda).
\end{equation}

We apply two methods to investigate upper and lower bounds for the extreme zero $x_{1,n}(\lambda)$
of  $C_{n}^{(\lambda)},$ $\lambda$ fixed, $-3/2 < \lambda < -1/2.$  One method
emanates from a suitably chosen mixed three-term recurrence relation and the other from the
Euler-Rayleigh technique discussed in \cite{im} where it is used to derive bounds
for the smallest real zero of a power series or  polynomial, with exclusively real zeros
in terms of the coefficients of the series or polynomial.

\noindent  Since the approach using mixed three term recurrence relations involves interlacing
properties of zeros of polynomials, we recall the definition:

\begin{Definition}
Let $\{p_n\}_{n=0}^\infty$ be a sequence of polynomials and suppose the zeros of $p_n$ are real
and simple for each $n \in {\mathbb N}.$ Denoting the zeros of $p_n$ in decreasing order by
$x_{n,n}<\dots <x_{2,n}<x_{1,n},$ the zeros of $p_n$ and $p_{n-1}$ are interlacing if,
for each $n \in {\mathbb N},$
\begin{equation}x_{n,n}<x_{n-1,n-1}<\dots <x_{2,n}<x_{1,n-1}<x_{1,n}.\label{1}
\end{equation}
\end{Definition}

The interlacing of the zeros of orthogonal polynomials of consecutive degree, $p_n$ and $p_{n-1}$  is a well known classical result
\cite[\S 3.3]{Sze}.
\medskip

\noindent Quasi-orthogonal polynomials arise in a natural way in the context of classical orthogonal polynomials that depend on one or more parameters. The concept of quasi-orthogonality of order $1$ of a sequence of polynomials was introduced by
Riesz \cite{Rie} in connection with the moment problem. Fej\'{e}r \cite{Fej}
considered quasi-orthogonality of order $2$ and the general case was studied by
Shohat \cite{Sho} and many other authors including Chihara \cite{Chi}, Dickinson \cite{Dic},
Draux  \cite{Dra} and Maroni \cite{ mar1, mar2, mar3}.

\medskip \noindent The definition of quasi-orthogonality of a sequence of polynomials is
the following:
\begin{Definition}
 Let $\{q_n\}_{n=0}^\infty$ be a sequence of polynomials with degree $q_n = n$ for each $n \in
 {\mathbb N}$. For a positive integer $r < n,$ the sequence $\{q_n\}_{n=0}^\infty$ is
 quasi-orthogonal of  order $r$ with respect to a positive Borel measure $\mu$ if
\begin{equation}
\int x^k q_n(x) \, d{\mu} (x) =0 \quad \mbox{for} \quad k=0,\ldots,n-1-r.
 \label{4} \end{equation}
\end{Definition}

\noindent If \eqref{4} holds for $r=0$, the sequence $\{q_n\}_{n=0}^\infty$ is orthogonal with
respect to the measure $\mu$.

\medskip \noindent
The sequence of
ultraspherical polynomials $\{C_{n}^{(\lambda)}\} _{n=0}^\infty,$ is orthogonal on $(-1,1)$ with
respect to the weight function $(1-x^2)^{\lambda -{1}/{2}}$ when $\lambda > -{1}/{2}.$
As $\lambda$ decreases below $-1/2,$ the zeros of $C_{n}^{(\lambda)}$ leave
the interval of orthogonality $(-1,1)$ in pairs through the endpoints as described in our opening paragraph. A full description of the departure of the zeros can be found in \cite{DrDu}.
For the range  $-3/2 < \lambda < -1/2,$ it is shown in  \cite[Theorem 6]{BDR} that the sequence
$\{C_{n}^{(\lambda)}\}_{n=0}^\infty$ is  quasi-orthogonal of order $2$  with respect to the
weight function $(1-x^2)^{\lambda +{1}/{2}};$ the polynomial
$C_{n}^{(\lambda)}$  has $n-2$ real, distinct zeros, its smallest zero is $<-1,$ its largest
zero is $>1$ and, for each $n \in{\mathbb N}$, the zeros of  $C_{n}^{(\lambda)}$ interlace with
the zeros of the (orthogonal) polynomial $C_{n-1}^{(\lambda+1)}.$

\medskip

\noindent In this paper, we derive upper and lower bounds for the
two (symmetric about the origin) zeros of $C_{n}^{(\lambda)}$ that lie outside $(-1,1)$ using two different methods.

\medskip

We use the notation $x_{k,n}(\lambda)$ for the $k$th zero, in {\it decreasing} order,
of $C_n^{(\lambda)}(x).$
It follows from \cite[(4.7.30)]{Sze} that $C_{n}^{(\lambda)}(x) \equiv 0$, for the $\lambda$-values
$\lambda^* = 0,-1,\dots,-\lfloor (n-1)/2\rfloor$. In this paper, we consider non-trivial zeros only. For each value of $\lambda^*$,
$C_n^{(\lambda^*)}(x)$  has $n$ non-trivial zeros defined  by
$$ x_{k,n}(\lambda^*) = \lim_{\lambda \to \lambda^*} x_{k,n}(\lambda).$$

\section{Bounds for the largest zero of $C_n^{(\lambda)}, \; -3/2 < \lambda < -1/2$.  An approach using
mixed three-term recurrence relations.}

\noindent The co-primality of two polynomials in a
sequence is an important factor when considering interlacing properties of their zeros.
Conversely, the inequalities satisfied by zeros of polynomials that are interlacing may
provide information about the co-primality of the polynomials. The interplay between
the co-primality of two polynomials and the interlacing of their (real) zeros facilitates
information about upper and lower bounds of the extreme zeros of
$C_n^{(\lambda)}, \; -3/2 < \lambda < -1/2$.
\begin{Theorem} \label{2.1}
Suppose that $\{C_n^{(\lambda)}\}_{n=0}^\infty$ is the sequence of ultraspherical polynomials
where $\lambda$ is fixed and lies in the range $-3/2 < \lambda < -1/2.$  For each
$n \in{\mathbb N},$ the largest zero, $x_{1,n},$ of $C_{n}^{(\lambda)}$, satisfies

\begin{equation} \label{3trru}
x_{1,n} <  \Big(  1+\frac{(2\lambda+1) }{(n-1) } \Big)^{-1/2}.
\end{equation}
\end{Theorem}
\noindent {\bf{Proof.}}
From \cite[Eqn.(16)]{Dr},
 we have
\begin{equation}
4 \lambda(\lambda+1){(1-x^2)}^2C_{n-2}^{(\lambda+2)}(x) =
(2\lambda+n)[x^2(n+2\lambda+1)-n)]C_n^{(\lambda)}(x)-(2\lambda+1)(n+1)x]
C_{n+1}^{(\lambda)}(x).\label{35}
\end{equation}
Evaluating \eqref{35} at the largest two zeros $x_{1,n+1} > 1 > x_{2,n+1}$ of
$C_{n+1}^{(\lambda)},$ we have
\begin{eqnarray}
 16{\lambda}^2(\lambda+1)^2(1-{x_{1,n+1}}^2)^2 (1-{x_{2,n+1}}^2)^2
 C_{n-2}^{(\lambda+2)}(x_{1,n+1})C_{n-2}^{(\lambda+2)}(x_{2,n+1})=
  \nonumber \\
(2\lambda+n)^2(2\lambda+n+1)^2C_n^{(\lambda)}(x_{1,n+1})C_n^{(\lambda)}(x_{2,n+1})
  \;\;\;\;\;\;\;\; \;\;\;\;\; \;\;\;\;\;\;\;\; \;\;\;\;\; \label{36} \\
\times \;\;[x_{1,n+1}^2-\frac{n}{n+2\lambda+1}] \nonumber
[x_{2,n+1}^2-\frac{n}{n+2\lambda+1}] .\end{eqnarray}
\noindent We know from \cite[Theorem 2.1(i)]{DrMu16}, with $n$ replaced  by $n+1$, that $C_{n}^{(\lambda)}$ does not change sign between
the two largest zeros of $C_{n+1}^{(\lambda)}.$ Further, since $C_n^{(\lambda)}$ and
$C_{n+1}^{(\lambda)}$ are co-prime, the only possible common zeros of
$C_{n-2}^{(\lambda+2)}(x)$ and $C_{n+1}^{(\lambda)}(x)$ are at the values
$x^2 = n/(n+2\lambda+1) = 1-(2\lambda+1)/(n+2\lambda+1)$ which have absolute value $>1$
for each $\lambda$ satisfying $-3/2 < \lambda < -1/2$ and $n \geq 3.$ Since, because of orthogonality, all the zeros of
$C_{n-2}^{(\lambda+2)}$ lie in $(-1,1)$ we deduce that  $C_{n-2}^{(\lambda+2)}$
and $C_{n+1}^{(\lambda)}$ are co-prime for each $n \in{\mathbb N},$ and each $\lambda$
satisfying
$-3/2 < \lambda < -1/2.$

\medskip
\noindent 
From \cite[Theorem 2.7(ii)]{DrMu16} with $n$ replaced by $n+1,$ we know that the zeros of
$(1-x^2) C_{n-2}^{(\lambda+2)}(x)$ and $C_{n+1}^{(\lambda)}(x)$ are interlacing.
Since the point $1$ lies between the two zeros $x_{1,n+1}$ and $x_{2,n+1},$
we see that in \eqref{36}, both the left-hand side and the product of the first four factors
on the right are positive and therefore
$\sqrt{{n}/(2\lambda +n+1)} \notin (x_{2,n+1},x_{1,n+1}).$ Also,
$\sqrt{{n}/ (2\lambda +n+1)}>1$ so that
$x_{1,n+1} < \sqrt{{n}/ (2\lambda +n+1)}.$ Finally, replacing $n$ by $n-1,$
we have the stated result.

\begin{Theorem} \label{2.2}
Suppose that $\{C_n^{(\lambda)}\}_{n=0}^\infty$ is the sequence of ultraspherical polynomials
where $\lambda$ is fixed and lies in the range $-3/2 < \lambda < -1/2.$ For each
$n \in{\mathbb N},$ the largest zero, $x_{1,n},$ of $C_{n}^{(\lambda)}$
satisfies
\begin{equation} \label{3trrl}
 x_{1,n} > \Big(  1+\frac{(2\lambda+1)(2\lambda+3) }{(n-1)(n+2\lambda+1) }    \Big)^{-1/2}.
 \end{equation}
\end{Theorem}

\noindent {\bf Proof.} From \cite[Eqn.(18)]{Dr}, we have
\begin{eqnarray}
8\lambda(\lambda+1)(\lambda+2){(1-x^2)}^3C_{n-2}^{(\lambda+3)}(x) =
\;\;\;\; \;\;\;\;\;\;\;\; \;\;\;\;\;\;\;\; \;\;\;\;\;\;\;\; \;\;\;\; \;\;\;\;\;\;\;\; \;\;\;\;\;\;\;\; \;\;\;\;\;\;\;\;\;\;\;\; \;\;\;\;  \label{37} \\
(2\lambda+n)[x^2[n(n+2\lambda+2)+(2\lambda+1)(2\lambda+3)]-n(n+2\lambda+2)]C_n^{(\lambda)}(x)
- g(x)C_{n+1}^{(\lambda)}(x), \nonumber
\end{eqnarray}
where $g(x)$ is a polynomial in $x$.
Evaluating \eqref{37} at the largest two zeros $x_{2,n+1}< 1 < x_{1,n+1}$ of $C_{n+1}^{(\lambda)},$
we know from \cite[Theorem 2.1]{DrMu16} that $C_{n}^{(\lambda)}$
does not change sign between the two largest zeros of $C_{n+1}^{(\lambda)}$
and, in addition,  $C_n^{(\lambda)}$ and $C_{n+1}^{(\lambda)}$ are co-prime. Therefore,
the only possible common zeros of
$C_{n-2}^{(\lambda+3)}(x)$ and $C_{n+1}^{(\lambda)}(x)$ are when
$$x^2 = \frac{n(n+2\lambda+2)}{n(n+2\lambda+2) +(2\lambda+1)(2\lambda+3)} = 1-\frac{(2\lambda+1)(2\lambda+3)}{n(n+2\lambda+2) +(2\lambda+1)(2\lambda+3)},$$ which is $>1$
for each $-3/2 < \lambda < -1/2$ and $n \geq 3.$ Since all the zeros of $C_{n-2}^{(\lambda+3)}$
lie in $(-1,1)$ we deduce that  $C_{n-2}^{(\lambda+3)}$
and $C_{n+1}^{(\lambda)}$ are co-prime. Also, the factor ${(1-x^2)}^3$ is positive at
$x = x_{2,n+1}$ and negative at $x= x_{1,n+1}.$ The proof now proceeds in the same way as the proof of
Theorem 2.1 by using an interlacing property  \cite[Theorem 2 (ii)(d)]{Dr} involving the zeros
of $C_{n-2}^{(\lambda+3)}$ and $C_{n+1}^{(\lambda)}$.

\section{The Euler-Rayleigh method}
					
A method described in \cite{im}  and applied to Laguerre and $q$-Laguerre polynomials
in \cite{gm}, may also be applied to prove bounds for the largest zero of $C_n^{(\lambda)},$
 $\lambda$ fixed, $-3/2 < \lambda < -1/2$. It may, in fact, also be applied to finding bounds for the largest zero of $C_n^{(\lambda)}$
 in the orthogonal cases when $\lambda > -1/2$.
\medskip

 The idea is as follows.  If  a polynomial $f(t) = \sum_{k=0}^na_nt^n$ has all its zeros real at the points
 \begin{equation}  t_1  < t_2 < \cdots < t_n, \label{1zn} \end{equation}
 and we use the notation
 $$S_j = \sum_{k=1}^n t_k^{-j},$$
 then the sums $S_j$ can be expressed in terms of the coefficients by
 $S_1 = -a_1, S_2 = -2a_2 +a_1^2,$ and, in general
$$ S_j = -na_n - \sum_{i=1}^{j-1}a_i S_{j-i}.$$
Two special cases can be considered.

\begin{Lemma}  \label{3.1}If $ 0 < t_1  < t_2 < \cdots < t_n,$ then $S_m >0, m=1,2,\dots,$ and
$$ S_m^{-1/m} < x_1 < S_m/S_{m+1},\;m=1,2,\dots$$
where the lower limits increase and the upper limits decrease with increasing $m$.
\end{Lemma}

This is \cite[Lemma 3.2]{im}.
\begin{Lemma} \label{3.2}
 If $ t_1 < 0 < t_2 < \cdots < t_n,$ and $a_1>0$, then $S_m <0$ for odd $m$,
\begin{equation}    -|S_{2m-1}|^{-1/(2m-1)} < t_1 <  -S_{2m}^{-1/(2m)}  <  S_{2m-1}/S_{2m}, \;m=1,2,\dots, \label{3.31}
\end{equation}
and \begin{equation}    S_{2m} /S_{2m+1} < t_1, \;m=1,2,\dots, \label{3.32}.
\end{equation}
\end{Lemma}

This is \cite[Lemma 3.3]{im} with the addition of the the remark at the
end of \cite[\S 3]{im}.

\medskip

We apply Lemma \ref{3.2} in the special case $m=1$ to the representation \cite[(4.7.6)]{Sze}
\begin{equation}
  C_n^{(\lambda)}(x) = \binom{n+2\lambda-1}{n}{} _2F_1\left(-n,n+2\lambda;\lambda +1/2;
t\right),\;\;t =  \frac{1-x}{2}. \label{uash}
\end{equation}
Since we are assuming that $-3/2 < \lambda < -1/2$ with $\lambda$ fixed, it follows from \cite[Cor. 2]{BDR}
that $x_n < -1 < x_{n-1} < \cdots < x_2 < 1 < x_1$ so the corresponding $t$-zeros satisfy
$t_1 < 0 < t_2 < \dots < t_n$.
With $f(t) = {} _2F_1\left(-n,n+2\lambda;\lambda +1/2;
t\right)$, we have
$$ S_1 = {\frac {n \left( n+2\,\lambda \right) }{\lambda+1/2}},
\;\;\;S_2= -S_1 \left[ \frac{2(n-1)(n+2\lambda+1)}{2\lambda+3}-S_1\right].$$

Applying the case $m=1$ of \eqref{3.31},  leads to
$$ -|S_1|^{-1} < t_1 < -S_2^{-1/2} < S_1/S_2.$$
  When converted to inequalities for $x_1 = 1-2t_1$, these inequalities give
$$
1-2S_1/S_2 < 1+2S_2^{-1/2} < x_{1,n}(\lambda) < 1 +2/|S_1, $$
or equivalently:
\begin{Theorem} \label{3.3} For each fixed $\lambda,$ $-3/2 < \lambda < -1/2$,
\begin{eqnarray}
\left[ 1+\frac{(2\lambda+1)(2\lambda+3)}{2(n-1)(n+2\lambda+1)}   \right]^{-1}
<
1-   \frac{2(2\lambda+1)\sqrt{2\lambda+3}}{\sqrt{ n ( 2\,\lambda+n ) ( 4\,{\lambda}^{2}+4\,n\lambda+2\,{n}^{2}
+4\,\lambda+1)}}  \nonumber\\
 < \;\;x_{1,n}(\lambda)
 < \;\;1-
  \frac{2\lambda+1} {n ( n+2\lambda) }. \;\;\;\; \;\;\;\; \;\;\;\; \;\;\;\; \;\;\;\; \;\;\;\; \label{1010}
\end{eqnarray}
\end{Theorem}

\section{Comparison of bounds}
It is of interest to compare the bounds in Sections 2 and 3, as well as to compare them with similar bounds valid for $\lambda > -1/2$.
\medskip

For $n \ge 6$, the upper bound in \eqref{1010} is sharper than that given by Theorem \ref{2.1}.  That is 
\begin{equation}
1- \frac{2\lambda+1} {n ( n+2\lambda) } <  \sqrt{\frac{(n-1)}{ (2\lambda +n)}},\;\; -\frac32 \le \lambda < -\frac12,\; n \ge 6, \label{comp1}
\end{equation}
with equality for $\lambda = -1/2$. The required inequality can be written $f(\lambda) > 0$ where %
$$ f(\lambda) = n^2(n+2\lambda)(n-1)-(n^2+2\lambda n-2\lambda -1)^2,$$
$$f'(\lambda) =2(n-1)(-n^2-4\lambda n+4\lambda +2), \;f''(\lambda) = -8(n-1)^2 <0.$$
Thus $f'(\lambda)$ is decreasing, and since $f'(-3/2) = -2(n-1)((n-3)^2-5) < 0$ we find that $f(\lambda)$ is decreasing on $(-3/2,-1/2)$. Finally since $f(-1/2)=0$, we see that $f(\lambda) >0$ for $-3/2 < \lambda < -1/2$.
\medskip
 
 For $n \ge 2$,   the lower bound given by Theorem \ref{2.2}
is sharper than the smaller lower bound in \eqref{1010}.  This statement is equivalent to the inequality
\begin{equation}
1+\frac{(2\lambda+1)(2\lambda+3)}{2(n-1)(n+2\lambda+1)}
<  \left[  1+\frac{(2\lambda+1)(2\lambda+3)}{(n-1)(n+2\lambda+1)}\right]^{-1/2}
\end{equation}
The number $a=[(2\lambda+1)(2\lambda+3)]/[(n-1)(n+2\lambda+1)]$ satisfies $-1 < a < 0$
for the values of $n$ and $\lambda$ concerned and it is a simple matter to show that, in this case, $1+a/2 <
1/\sqrt{1+a}$.

\medskip
We can use Lemma \ref{3.1} to show that the inequalities (Theorem \ref{3.3})
\begin{equation}
\left[ 1+\frac{(2\lambda+1)(2\lambda+3)}{2(n-1)(n+2\lambda+1)}   \right]^{-1}
 < \;\;x_{1,n}(\lambda)
 < \;\;1-
  \frac{2\lambda+1} {n ( n+2\lambda) } \label{uandl}
\end{equation}
continue to hold for $\lambda > -1/2$. with equality  for $\lambda = -1/2$.

\section{A continuity-based proof of interlacing}
In  \cite[Theorem 2.1]{DrMu16}, we showed that for $-3/2 < \lambda < -1/2$,
the zeros of $(1-x^2)  C_{n-1}^{(\lambda)}(x)$ interlace with the zeros
of $C_{n}^{(\lambda)}(x)$.  The same result holds for $\lambda > -1/2$ as
can be seen by adding the two points $\pm 1$ to the well-known interlacing
of the zeros (all in $(-1,1)$) of $C_{n-1}^{(\lambda)}(x)$ and $C_{n}^{(\lambda)}(x)$.
Here we show how the interlacing for $\lambda > -1/2$ may be used to get the interlacing
for $-3/2 < \lambda >-1/2$.

\medskip
The proof is based on the idea that interlacing between the zeros of $C_n^{(\lambda)}(x)$ and $C_{n+1}^{(\lambda)}(x)$ can break down or change only when these two functions have a common zero.
First of all, we prove some lemmas.

\begin{Lemma}
The only common zeros of $C_n^{(\lambda)}(x)$ and $C_{n+1}^{(\lambda)}(x)$ occur where
$|x| =1$ and $\lambda = -1/2, -3/2,\dots, -\lfloor n/2 \rfloor$.
\end{Lemma}

{\bf Proof.}
Because of parity considerations, the polynomials cannot have a common zero at $0$. Let $\lambda = \lambda_0$ be fixed.
The recurrence relation \cite[(4.7.17)]{Sze}
\begin{equation}
nC_{n}^{(\lambda)}(x) = 2(n+\lambda-1) xC_{n-1}^{(\lambda)}(x) - (n+2\lambda-2)C_{n-2}^{(\lambda)}(x),
\;n=2,3,4,\dots
\end{equation}
show that, if $C_n^{(\lambda_0)}(x)$ and $C_{n-1}^{(\lambda_0)}(x)$ had a common zero $x_0$ for a value of
$x$ satisfying $0 < |x| \ne1$, then $x_0$ would also be a zero of
 $C_{n-2}^{(\lambda_0)}(x)$.
 Repeating this argument we would find that $C_1^{(\lambda_0)}(x_0) =0,$
 which is impossible.
\medskip

Thus the only possible common zeros of $C_n^{(\lambda)}(x)$ and $C_{n+1}^{(\lambda)}(x)$
occur for $x= \pm 1$.   From \cite[Ch. 5]{hartman} the zeros of 
$C_n^{(\lambda)}$ are continuous functions of $\lambda$
and, as $\lambda$ varies, the only values of $\lambda$ for which interlacing breaks down are where  $C_n^{(\lambda)}(x)$ and $C_{n+1}^{(\lambda)}(x)$ have common zeros.
  Now \cite[(4.7.6)]{Sze}, $C_n^{(\lambda)}(x)$ is a multiple of
${}_2F_1(-n, n+2\lambda; \lambda + \frac12; \frac{1-x}{2}) $, it vanishes at $1$ only for
$\lambda = -1/2, -3/2,\dots, -\lfloor n/2 \rfloor$ and $C_{n+1}^{(\lambda)}(x)$ vanishes at $1$ only for
$\lambda = -1/2, -3/2,\dots, -\lfloor (n+1)/2 \rfloor$.  The statement of the Lemma follows.
\endproof

\begin{figure} \label{n=8,9}
\centering
\includegraphics[width=4.5in, height=3.5in,angle=0]{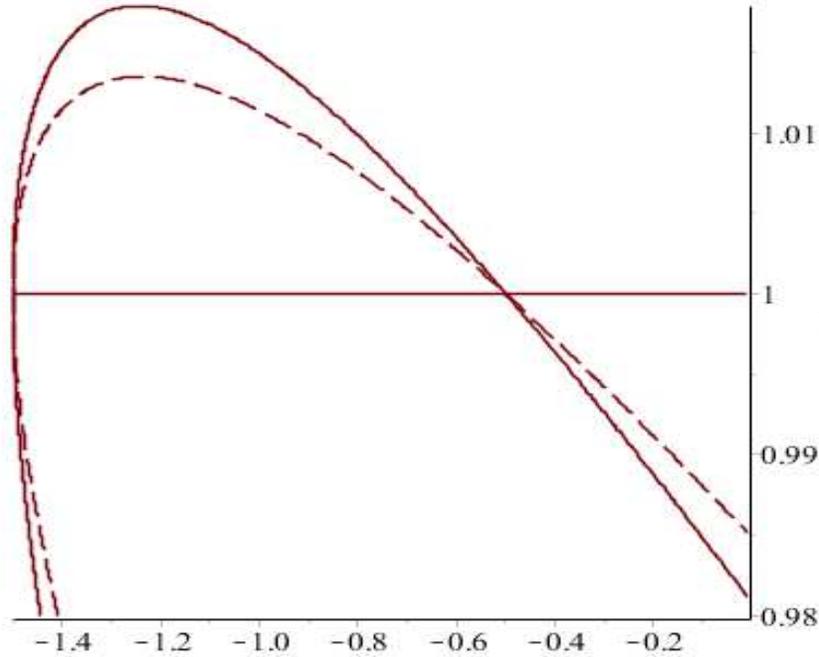}
\caption{Largest zeros of $C_8^{(\lambda)}(x)$ (solid curves) and
$C_{9}^{(\lambda)}(x)$ (dashed curves), as functions of $\lambda, \;-3/2 < \lambda < 0.$}
 \end{figure}

 \begin{Lemma} \label{5.2}
 If $x_{1,n}(\lambda)$ is the largest zero of $C_{n}^{(\lambda)}(x)$, then
 $$  x_{1,n}'(-1/2) = -2/(n^2-n). $$
 In particular, $x_{1,n}'(-1/2)$ increases with $n,\;n=2,3,\dots$.
 \end{Lemma}

 {\bf Proof.} The functions giving upper and lower bounds in \eqref{uandl} both have derivatives
 equal to $ -2/(n^2-n)$ at the point $\lambda = -1/2$.  Then \eqref{uandl} implies that $x_{1,n}'(-1/2)$
 must have the same value.
 \endproof

We are now ready to show how to extend the property of the interlacing of real zeros of
$(1-x^2)C_{n}^{(\lambda)}(x)$ and $C_{n+1}^{(\lambda)}(x)$ from the case $\lambda > -1/2$ to the case $-3/2< \lambda < -1/2$.  To keep things simple we consider positive zeros only.

We write the zeros in decreasing order.  For $\lambda > -1/2$, we have
\begin{equation}
1 >  x_{1,n}(\lambda) >  x_{2,n}(\lambda) >  x_{3,n}(\lambda) > \dots,
\end{equation}
and the interlacing property
\begin{equation}\label{22}
1 >  x_{1,n+1}(\lambda) >  x_{1,n}(\lambda) >  x_{2, n+1}(\lambda) >  x_{2, n}(\lambda) > \dots.
\end{equation}
 From Lemma 5.2,
 the slope of  $x_{1,n}(\lambda)$, for $\lambda = -1/2$, is a negative increasing function of
$n$ and so the functions $1,  x_{1, n+1}(\lambda) $ and $x_{1,n}(\lambda)$ have their order reversed in the inequality \eqref{22} as
$\lambda$ passes through the value $-1/2.$ 
Thus, for $-3/2 < \lambda < -1/2$, we get
\begin{equation}\label{23}
\boxed{ x_{1,n}(\lambda)  >  x_{1,n+1}(\lambda) > 1} >  x_{2,n+1}(\lambda) >  x_{2,n}(\lambda) > \dots,
\end{equation}
with the other inequalities remaining the same, since the values $\lambda =-1/2$ and $x=1$ constitute the only double zero within the range
considered.

\medskip

Figure 1,  produced using Maple, illustrates the change of order of  $x_{1,n}(\lambda), x_{1,n+1}(\lambda) $ and $1$ as lambda passes through $-1/2$ in the case $n=8$.

\medskip

{\bf Acknowledgments.}  Kathy Driver's research was supported by the National Research Foundation of South Africa.

\medskip

{Department of Mathematics and Applied Mathematics, University of Cape Town, Private Bag X1, Rondebosch 7701, South Africa.   email:  Kathy.Driver@uct.ac.za}
 
 \medskip
 
 {Department of Mathematics and Statistics, York University, Toronto, ON   
   M3J 1P3, Canada.  email: muldoon@yorku.ca}

\end{document}